\documentclass{amsart}
\usepackage{float}
\usepackage{amsfonts,amssymb,amsmath,url,amsthm}

\newtheorem{theorem}{Theorem}%[section]
\newtheorem{lemma}{Lemma}
\newtheorem{corollary}{Corollary}
\newtheorem{conjecture}{Conjecture}

\newtheorem*{theorem*}{{\bf{Theorem}}}
\newtheorem*{corollary*}{{\bf{Corollary}}}

\theoremstyle{remark}

\newtheorem{remark}{Remark}
\newtheorem*{example*}{{\bf{Example}}}
\newtheorem*{lemma*}{{\bf{Lemma}}}
\newtheorem*{remark*}{{\bf{Remark}}}

\newcommand{\Z}{{\mathbb{Z}}}

\newcommand{\F}{{\mathbb{F}}}
\newcommand{\Q}{{\mathbb{Q}}}

\newcommand{\ord}{{\rm{ord }}}

\title[Almost all primes have a multiple of small Hamming weight 
]{Almost all primes have a multiple of small Hamming weight}
\author{
Christian Elsholtz}\thanks{
Institute of Analysis and Number Theory,
Graz University of Technology
Kopernikusgasse 24,
A-8010 Graz, Austria.
{\it elsholtz@math.tugraz.at}}

\begin{document}
\maketitle

\begin{abstract}
Recent results of Bourgain and Shparlinski imply
that for almost all primes $p$ there is a multiple $mp$ 
that can be written in binary as 
\[mp= 1+2^{m_1}+  \cdots +2^{m_k}, \quad 1\leq m_1 < \cdots < m_k,\]
with $k=66$ or $k=16$, respectively. We show that $k=6$ (corresponding to
Hamming weight $7$) suffices.

We also prove there are infinitely many  primes $p$ 
with a multiplicative subgroup $A=<g>\subset \F_p^*$, for some
 $g \in \{2,3,5\}$, of size $|A|\gg p/(\log p)^3$, where the sum-product set
$A\cdot A+ A\cdot A$ does not cover $\F_p$ completely.
\end{abstract}

\section{Introduction}
Recently Shparlinski \cite{Shparlinski:2003} initiated the
study of prime divisors of ``sparse integers'', i.e. of integers
which only use a few non-zero digits in a $g$-ary representation.
Let $g \geq 2$ and $k\geq 1$ be integers, and let ${\mathcal D}=(d_i)_{i=0}^k$
be a sequence of $k+1$ nonzero integers. Let 
${\mathcal S}_{g,k}({\mathcal D})$ be the
set of integers $n$ of the form
\begin{equation}{\label{eq:mainform}}
n=d_0+d_1g^{m_1}+ \cdots +d_k g^{m_k},
\end{equation}
 where $1\leq m_1 < \cdots < m_k$. Here the $d_i$ can be thought of as 
(but are not necessarily) the
 $g$-ary digits $d_i \in \{1, \ldots, g-1\}$ of $n$.
Using exponential sums Shparlinski \cite{Shparlinski:2003} proved: 
for any fixed $\delta\in (0,
1/2)$ and every $k> \max\{15, \frac{1}{\delta}-1\}$, 
for sufficiently large $X$, for almost all primes $p \leq X$,
(i.e. for all but $o(\pi(X))$ primes, where $\pi(X)$ denotes the number of
primes $p \leq X$), there exists $n \in {\mathcal S}_{g,k}({\mathcal D})$ with
$\log n \ll X^{1/2+\delta}$, such that $p \mid n$ holds.

Using new bounds of very short exponential sums
Bourgain \cite{Bourgain:2005GAFA, Bourgain:2004} managed to prove the existence
of much smaller $n$: Let $\delta >0$, and $k>k_0(\delta)$. Let $X$ 
be sufficently large, then for all but $o(\pi(X))$ primes $p \leq X$ 
there exists  $n \in {\mathcal S}_{g,k}({\mathcal D})$ with
$\log n \ll X^{\delta}$, such that $p \mid n$ holds.

Shparlinski \cite{Shparlinski:2008}
writes, working out $k_0$ from Bourgain's approach would
take ``significant efforts'', and with an alternative approach 
using incomplete exponential sums ``on average over primes''
he worked out an explicit admissible value:
Let $0 < \delta < 1/2$, then one can take 
\[k_0(\delta)= \frac{16}{\delta^2} +1.\]

One may wonder, how much this can be improved.
It is clear that $k \geq \frac{1}{\delta}$ must hold:
The number of sparse integers of the type (\ref{eq:mainform})
is of order of magnitude $C_{k,g} (\log n)^k\ll X^{\delta k}$, and this
expression must be at least 
of the same order of magnitude as the number of primes, $\frac{X}{\log X}$.

It seems there is no result in the literature that tries to minimize $k$,
even when allowing a larger value $n$. It follows from the results above that
$k=16$, respectively $k=66$ are admissable.
Here we will %concentrate on the most important case $g=2$ and 
show that actually $k=6$ suffices. 
As it turns out the proof combines essentially two
types of tools from the literature:\\
a) recent results on sumsets modulo a fixed prime $p$ 
from additive combinatorics,\\
b) an old result of Erd\H{o}s stating that, for most primes,
the powers of 2 generate at least $p^{1/2+o(1)}$ distinct residue classes, and
variants, e.g. due to
 Pappalardi \cite{Pappalardi:1996}, Erd\H{o}s and Ram Murty
\cite{ErdosandRamMurty:1999}.

Using these tools the proof is actually quite short. 
However it seems the new type of
view and the method of proof lead to considerable progress on the problem
mentioned.

\section{Discussion of the  case $g=2$.}

The Hamming weight 
of an integer counts the number of ones in binary representation.
During the last years
there have been a number of related investigations on different 
number theoretic questions, such as
primes $p$ with all positive 
values of $p-2^i$ composite (see Tao \cite{Tao:2011})
or quadratic nonresidues with small Hamming weight,
(see Dietmann, Elsholtz and Shparlinski \cite{DietmannElsholtzandShparlinski:2012}).

In this article we show that almost all primes have a nonzero
multiple with a bounded Hamming weight of at most $7$ 
(corresponding to $k=6$ above).

Let $T_k$ denote the set of those primes $p$ which divide some integer of
the form $2^{a_1}+\cdots + 2^{a_k}+1$, 
where $a_1, \ldots , a_k \in {\mathbb N}$, and let $T_k(X)$ denote the number
of such primes $p \leq X$.

A result of Hasse \cite{Hasse:1966}
says that the Dirichlet density of primes $p \leq x$ 
dividing a number of the form 
$2^a+1$, where $a\in {\mathbb N}$, is
$\frac{17}{24}$.  
This was later conjectured for natural density by Krishnamurty, and eventually
proved by Odoni \cite{Odoni:1981}.
In other words,
$T_1(X)\sim \frac{17}{24} \pi(X)\sim \frac{17}{24} \frac{X}{\log X}$.

There are extensions of this result to composite moduli, or to more general
sequences such as $a^k+b^k$
 due to Ballot \cite{Ballot:1995}, Moree \cite{Moree:1997},
Odoni \cite{Odoni:1981} and Wiertelak \cite{Wiertelak:1984}.

Ska\l ba \cite{Skalba:2004} conjectured, 
(see also Moree \cite{Moree:2012}):
\begin{conjecture}[Ska\l ba]{\label{con:skalba}}
$T_2(X)\sim\frac{X}{\log X}$.
\end{conjecture}

Ska\l ba proved a partial result towards this conjecture. 
To formulate it we introduce
some background first.
Let $\ord_p(a)$ denote the multiplicative order of $a$ in $\F_p^*$. 
Erd\H{o}s \cite{Erdos:1976}
conjectured that for each $0<c<1$ for almost all primes
${\ord}_p(2)>p^c$ holds, and proved this for $c\leq \frac{1}{2}-o(1)$.
There have been several attempts to address Erd\H{o}s's conjecture.
 Indeed, on the Generalized Riemann
Hypothesis it is known that for almost all primes $p$:
${\rm ord}_p(2)>\frac{p}{f(p)}$, where
$f$ tends (arbitrarily slowly) to infinity, 
see Pappalardi \cite{Pappalardi:1996}.
Unconditionally, refining work of Pappalardi \cite{Pappalardi:1996},
Erd\H{o}s and Murty \cite{ErdosandRamMurty:1999} proved that 
for almost all primes $p \leq x$ and for any function
$\varepsilon(x)$ tending to $0$ as $x$ tends to infinity,
\[ \ord_p(2) \geq p^{\frac{1}{2}+\varepsilon(p)} \]
holds. 

Ska\l ba gave a conditional proof of his conjecture:
\begin{theorem*}[Ska\l ba \cite{Skalba:2004}]
If ${\rm ord}_p(2)>p^{0.8}$, then there exist integers $a,b$
such that
\[2^a+2^b+1\equiv 0 \bmod p\]
holds. In particular, if
Erd\H{o}s's conjecture holds with $c=0.8$, 
then Conjecture \ref{con:skalba} holds.
\end{theorem*}
The proof made use of Weil bounds \cite{Weil:1949}
on the number of solutions on congruences.

As progress on Erd\H{o}s's conjecture 
has been very slow
it seems fair to say that Erd\H{o}s's conjecture with
$c=0.8$ is currently very far from a solution.
The currently known exponent $c=\frac{1}{2}$,
combined with Ska\l ba's method based on Weil estimates
seems insufficient to establish
$T_k(X)\sim\frac{X}{\log X}$, for any fixed $k$.

As mentioned in the introduction, 
the results by Shparlinski give $k=66$ and $k=16$, respectively.
We take a different approach to the problem.
Combining methods from additive and multiplicative number theory
we prove a result towards Skalba's conjecture, where
$2^a+2^b+1$ is replaced by $2^{a_1} + \cdots + 2^{a_6}+1$.
In other words, we prove that $T_6(x) \sim \frac{x}{\log x}$. This means 
that almost all primes $p$ have a nonzero multiple $mp$
with Hamming weight at most $7$.
Observe that
the number of integers $n \leq N$ with Hamming weight 
at most $k$ is about $\sum_{i=0}^k \binom{t}{i}$, where 
$t\sim \frac{\log N}{\log 2}$. 
Hence there are about $O((\log N)^6)$ odd integers
$n \leq N$ of Hamming weight at most $7$. 
For a given prime $p$ the corresponding multiple $mp$
of small Hamming weight can of course be much larger than the prime.
As the proof shows, the exponents $a_i$ are bounded above by $p-1$, and hence 
$mp\leq k2^{p-1}p$.
It seems possible to reduce the exponent by a small factor, 
following the techniques in
Cilleruelo and Zumalac\'{a}rregui \cite{CilleruleoandZumalacarregui}
and Garaev and Kueh \cite{GaraevandKueh:2002}, but we 
do not make an attempt to do this.

One may wonder if this result, 
that almost all primes have a very sparse 
binary multiple, has any practical application. 
Of course, addition and multiplication with sparse
integers should be faster than with arbitrary integers, 
but it seems one would need to find this sparse multiple in the first place
and there could be several more issues to overcome.

We also prove similar results, where
$2^a+2^b+1$ is replaced e.g.~by 
$2^{a_1} 3^{b_1} 5^{c_1} + 2^{a_2} 3^{b_2} 5^{c_2}+1$. A more general
formulation of these results is in the next section.

In the opposite direction Ska\l ba proved:
\begin{theorem*}[Ska\l ba \cite{Skalba:2004}]
Let $\Omega(w)$ denote the number of prime factors of $w$, with multiplicity.
If $\Omega(2^n-1)< \frac{\log n}{\log 3}$, then there exists a prime divisor
$p$ of $2^n-1$ such that for no pair of integers 
$(a,b): 2^a+2^b+1 \equiv 0 \bmod p$.
\end{theorem*}

We will extend this to a similar estimate with $2^a+2^b$ replaced by a
$k$-fold sum.
One may conjecture that for arbitrarily large
$k$ there are infinitely many primes $m$
such that no multiple of $m$ can be written as 
$\sum_{i=1}^k 2^{s_i}$ for
nonnegative integers $s_i$.
 This would of course follow if there exist infinitely many
Mersenne primes, but the condition
$\Omega(2^n-1)< \frac{\log n}{\log k}$
in Theorem {\ref{thm:Skalba2-extension}}
is much more modest.

In the last section, we study a restriction to sum product estimates.

\section{Results}
\subsection{Representation of residue classes, 
and approximations to Ska\l ba's conjecture
}
\begin{theorem}{\label{thm:6fold}}
Let $r\geq 2$ and $m\neq 0$ be fixed integers.
For almost all primes $p\leq x$ there is a solution of
\[r^{a_1}  + r^{a_2} + \cdots +r^{a_6}\equiv m \bmod p,\]
with integers $0 \leq a_1, \ldots, a_6 < p-1$.
More precisely, with $\varepsilon>0$,
the number of primes $p \leq X$ with no solution of
\[r^{a_1}  + r^{a_2} + \cdots +r^{a_6}\equiv m \bmod p\]
is 
$O_{r}(X^{22/23+\varepsilon})$.
\end{theorem}
With $r=2, m=-1$ we get the following approximation
to Ska\l ba's conjecture.
\begin{corollary}{\label{cor-sixpowers}}
The following holds: $T_6(X)\sim\frac{X}{\log X}$.
\end{corollary}
\begin{remark}
In the special case that there exists some $i$ with $2^i\equiv -1 \bmod p$
one can reduce the ``6'' in Theorem {\ref{thm:6fold}} above
(with a different size of the exceptional set)
to ``5'',
in view of recent work of Shkredov (Theorem 2 of \cite{Shkredov:2013}). 
\end{remark}
\begin{theorem}{\label{thm:hart-2fold}}
Let $r\geq 2$.
For almost all primes $p\leq X$ 
there are at least $p^{4/5}$ classes $m$ which can be represented by
\[r^{a_1}  + r^{a_2} \equiv m \bmod p,\]
with integers $0 \leq a_1, a_2 < p-1$.
\end{theorem}

Other approximations to Ska\l ba's conjecture are as follows:
\begin{theorem}{\label{thm:3fold}}
Let $r,s\geq 2$ be coprime integers
and let $m\neq 0 $ be a fixed integer.
For almost all primes $p\leq X$ there is a solution of
\[r^{a_1} s^{a_2} + r^{a_3} s^{a_4}+r^{a_5}s^{a_6}\equiv m \bmod p,\]
with integers $0 \leq a_1, \ldots, a_6 < p-1$.
\end{theorem}
\begin{theorem}{\label{thm:2fold}}
Let $r,s,t\geq 2$ be mutually coprime integers 
and let $m\neq 0$ be a fixed integer.
For almost all primes $p\leq X$ there is a solution of
\[r^{a_1} s^{a_2} t^{a_3} + r^{a_4} s^{a_5} t^{a_6}\equiv m \bmod p,\]
with integers $0 \leq a_1, \ldots, a_6< p-1$.
\end{theorem}
For completeness we also state the following,
which is a direct consequence of a result of
Schoen and Shkredov \cite{SchoenandShkredov:2011}.
\begin{corollary}{\label{thm:2foldconditional}}
Let $r\geq 2$.
If ${\rm ord}_p(r)> p^{3/4}$, then for some integers $a_1, a_2$
\[r^{a_1}+r^{a_2}+1\equiv 0 \bmod p.\]

\end{corollary}
In particular, if Erd\H{o}s's conjecture holds with $c=0.75$, 
then Conjecture \ref{con:skalba} holds.
This last result replaces $\frac{4}{5}$ in Ska\l ba's Theorem 
(mentioned above) by $\frac{3}{4}$.

\begin{remark}
Another approach to this problem could be the following one:
a) Baker and Harman \cite{BakerandHarman:1998} proved that 
for a positive proportion of primes,  
the largest prime factor $P$ of $p-1$ is large: $P(p-1)\gg p^{0.677}$.
This result is ineffective due to the use of Siegel zeros.\\
b) Harman \cite{Harman:2005} also gave an effective result 
(i.e.~not depending on a Siegel zero),
namely that for some positive (computable) 
density of the primes this holds with exponent 0.6105. 

From case a) one can deduce
%(see {\cite[Lemma 20]{KurlbergandPomerance:2005}}):
that for a positive  (ineffective) proportion of primes 
$\ord_p(2)\gg p^{0.677}$ holds. In view of the large prime factor 
we first consider the case in which $2$ 
generates only $O(p^{1/3})$ residue classes: here we know
from the Erd\H{o}s conjecture with $c=\frac{1}{2}$,
see also Lemma \ref{lem:erdos-murty}, 
that this case rarely 
happens. In the other case $2$ generates a subgroup of size at least 
$cp^{0.677}$ residue classes. See also 
{\cite[Lemma 20]{KurlbergandPomerance:2005}}.

This would imply (e.g.~by Lemma \ref{lemma:Schoenshkredov3})
that for these primes $2^a+2^b+2^c+1 \equiv 0 \bmod p$ has a solution.
In its current form this consequence would  
be much weaker than Hasse's result, 
namely that for a density of $17/24$ of the primes even   
$2^a+1 \equiv 0 \bmod p$ has a solution. Therefore 
we do not pursue this path further.

\end{remark}

\subsection{Primes without multiples of small Hamming weight}

As an extension of Ska\l ba's second theorem with regard to multiple sums,
and thus a limitation or partial converse to the type of results above,
we prove the following:
\begin{theorem}{\label{thm:Skalba2-extension}}
If
\begin{equation}%{\label{eq:mersennefactors}}
\Omega(2^n-1)< \frac{\log n}{\log k},
\end{equation}
 then there is some prime factor $q\mid 2^n-1$ such that
\[2^{a_1}+\cdots + 2^{a_{k-1}}+1 \not\equiv 0 \bmod q,\]
for all integer choices of $a_i$.
Hence all multiples of $q$ have Hamming weight at least $k+1$.
\end{theorem}

\subsection{A restriction on sum product estimates}

Suppose a $k$-fold sum $r^{a_1}+ \cdots + r^{a_k}$ covers all residue
classes, in particular the class $0$; 
here without loss of generality 
$a:=a_k\leq \cdots \leq a_1$. Dividing by $r^a$
one obtains that $r^{a_1-a}+ \cdots + r^{a_{k-1}-a}+1 \equiv 0 \bmod p$.

Also, the related question has been studied, about the minimal 
number $d$ such that the mixed $d$-fold sum-product set 
$d AB:= \{ a_1 b_1+ \cdots + a_d b_d: a_i \in A, b_i \in B\}$ covers all residue classes.

For random sets $A\subset \Z/p\Z$
with $|A| \geq C_{\varepsilon} p^{\frac{1}{2}+\varepsilon}$,
one may expect that $A+A=\Z/p\Z$.
See for example some discussion (with too optimistic conjectures) in
Hart and Iosevich \cite{HartandIosevich:2008},
Rudnev \cite{Rudnev} and 
Chapman, Erdo\u{g}an, Hart, Iosevich and Koh 
\cite{ChapmanErdoganHartIosevichandKoh:2012}.
The last authors observed that some kind of restriction must occur, they write
about subsets $A\subset \F_q$:
``Due to the misbehavior of the zero element it is not possible for 
$A \cdot A+ A \cdot A = \F_q$ unless
$A$ is a positive proportion of the elements of $\F_q$.''

It may be worth recalling that for a prime $p\equiv 3 \bmod 4$ the set 
$A$ of (nonzero) quadratic residues has positive density but 
$0 \not\in AA+AA$. The set of squares are of course an explicit
example of a large multiplicative subgroup of $\F_p^*$ 
with restricted sumsets in $\F_p$. On the other hand
this example of size $|A|=\frac{p-1}{2}$ is extremal,
as for any set $A$ with  $|A|\geq \frac{p+1}{2}$
the Cauchy-Davenport theorem guarantees that
$A+A=\F_p$.

The set of squares is generated by the square $g^2$ of a primitive root $g$.
But, as $g^2$ is not a fixed element,
following the view point of this paper one can ask if there exists such 
examples generated by a fixed element, e.g. $g=2$ and 
$A=\{2^i\bmod p: 1\leq i \leq \ord_p(2)\}$.
In short, the answer is: ``probably yes, 
but unconditionally we can only prove a slightly weaker result''.

Let us first describe an explicit and ``large'' example, where however we
cannot prove this type of example will occur infinitely often.
We then describe a slightly weaker example, which occurs quite frequently.

Let $a$ be fixed, and $p$ be a prime with
$\ord_p(a)=\frac{p-1}{2}$, where $\frac{p-1}{2}$ is odd.
%In particular $\left( \frac{a}{p}\right)=1$.
Let $A=B=\{a^i: 1\leq i \leq \frac{p-1}{2}\}$.
Then $|A||B| \sim \frac{p^2}{4}$ and $2AB=A^2+A^2= A+A$.
But $0\not\in A+A$, as otherwise $a^i\equiv -a^j\bmod p$. But $-1$ is not in
the multiplicative subgroup generated by $a$, as $a$ has odd order.
The size of $|A|=\frac{p-1}{2}$ 
is of course much larger than $p^{\frac{1}{2}+\varepsilon}$.
For $a=2$ such primes 
are $p=7,23,47,71,79,103,167,191,199,\ldots$.
Alternatively, one can take
 primes $p \equiv 3 \bmod 4$, with $2$ as primitive root,
and $A=\{4^i: 1\leq i \leq \frac{p-1}{2}\}$.
Such primes are
$3, 11, 19, 59, 67, 83, 107, 131, 139, 163, 179, \ldots$.
Artin's conjecture on primitive roots (known only on GRH, see Hooley
\cite{Hooley:1967})
predicts that such sets of primes 
have a positive density in the set of primes.
%The question if there are infinitely many such primes or not is in both cases
%closely related to Artin's conjecture on primitive roots.

Unconditionally,
we are able to prove the existence of an infinite number of primes
with a restricted sumset for the following slightly weaker variant:

\begin{theorem}{\label{thm:rudnevcounterexample}}
For a given prime $p$ let $s=s(p)$ and $w=w(p)$ be defined by 
$p=2^{s} w+1$, where $w$ is odd. 
For almost all positive square-free numbers $a>1$,
with at most three exceptions,
(or almost all primes, with at most 2 exceptions),
there exist positive constants $c$ and $c_1$,
and there exist at least $ c_1 \frac{x}{(\log x)^2}$ many primes $p \leq x$
such that the  multiplicative subgroup 
\[A=\{\left(a^{2^s}\right)^i: 1 \leq i \leq \ord_p(a^{2^s})\}\subset \F_p^*\]
has the following properties:

\begin{enumerate}
\item $|A|> c p/(\log p)^3$ and
\item $A\cdot A+ A \cdot A =A +A \neq \F_p$.
\end{enumerate}
\end{theorem}

Remark: in this situation $A=A \cdots A$, but as $|A \cdot A|$ is 
typically much larger than $|A|$ we state it in the form above.

Let us mention a related result by Alon and Bourgain 
\cite{AlonandBourgain:2014}:
there is an absolute constant $c > 0$ so that there are infinitely
many primes $p$ and a multiplicative subgroup $A \subset \F_p^*$
with $|A| \geq c p^{1/3}$ such that:  there are no $x, y,z \in  A$ 
with $x+y=z$.
 
\section{Proofs}

We observe that Erd\H{o}s's conjecture only makes a statement
 about the number of residue classes of the form $r^i\bmod p$, 
(for fixed $r$ and $p$),
but not about its
algebraic structure. 
Observing that these classes are a multiplicative subgroup
of $\F_p^{*}$ we now study sumsets of multiplicative subgroups 
in $\F_p^{*}$. In recent years there has been a considerable number of
papers concerned with sumsets modulo primes. Of the many results available
 we choose
those results that seem best suited for our application.

\subsection{Lemmas}
For the proof we essentially use an ingredient
(the first two lemmas below) coming from
multiplicative number theory and several very recent quantitative versions
of results of additive combinatorics.

\begin{lemma}[Erd\H{o}s-Murty 
{\cite[Theorem 5]{ErdosandRamMurty:1999}}\label{lem:erdos-murty}, 
Pappalardi \cite{Pappalardi:1996}]{\label{lem:erdosmurty}}
Let $\Gamma \subseteq \Q^*$ be a multiplicative subgroup of rank $d$. Suppose
that $\Gamma$ is generated by the mutually coprime numbers $b_1, \ldots , b_d$.
For all primes $p$ not dividing the denominators of $b_1, \ldots , b_d$ we
define $f_{\Gamma}(p)$ to be the order of $\Gamma \bmod p$.
Let $\varepsilon(x)$ be any function tending to zero as $x \rightarrow \infty$.
For all but $o(\frac{x}{\log x})$ primes $p \leq x$:
\[ f_{\Gamma}(p) \geq p^{\frac{d}{d+1}+\varepsilon(p)}.\]
\end{lemma}

The number of exceptional primes can be bounded 
from above in a more precise way;
see for example Indlekofer and Timofeev 
\cite{IndlekoferandTimofeev:2002} and Ford \cite{Ford:2008}.
However for our situation, a method by Matthews, worked out in detail
by Pappalardi, is suitable.

Here $b_1, \ldots , b_d$, are multiplicatively 
independent integers, not being squares, and not being $\pm 1$.

\begin{lemma}[Pappalardi {\cite[Lemma 1.2]{Pappalardi:1996}}, 
Matthews \cite{Matthews:1982}]{\label{lem:Pappalardi}}
Suppose that $d$ is a function of $t$ such that $dt^{-d}$ is bounded.
Then
\[ \# \{ p: f_{\Gamma}(p) \leq t\}
\ll \frac{t^{1+1/d}}{\log t} 2^d d \sum_{i=1}^d \log b_i\]
uniformly with respect to $t, d$ and $\{b_1, \ldots, b_d\}$.
\end{lemma}
Our application will only need the special case $d =1$.

A slightly weaker version of Lemma \ref{lem:erdos-murty} goes back 
to Erd\H{o}s \cite{Erdos:1976}, and one can give the following
very simple proof:

The sequence $(2^i \bmod p), i=1, \ldots , n$ is periodic.
The order $\ord_p(2)$ is by definition the length of the period,
and one has $2^{\ord_p(2)}\equiv 1 \bmod p$.
This implies that
\[ \prod_{p< y} p^{\lfloor \frac{n}{\ord_p(2)} \rfloor} 
\leq \prod_{i=1}^n (2^i-1) \leq 2^{n(n+1)/2},\]
where the product on the left hand side runs over the primes, up to some level
$y$ specified below.
Taking logarithms one sees that
\begin{equation}{\label{product-argument}}
\sum_{p< y} 
\log p \lfloor \frac{n}{\ord_p(2)} 
\rfloor \leq C n^2,
\end{equation}
for some positive constant $C$. Assuming that 
$\ord_p(2)\leq y^{\frac{1}{2}-\varepsilon}$ 
for at least $\delta \frac{y}{\log y}$
of the primes $p\leq y$, (for some $\delta >0$),
then
\[ \sum_{p< y} \log p \left( \frac{n}{\ord_p(2)} -1\right)
\geq  \frac{\delta y}{\log y} \left( \frac{n}{y^{\frac{1}{2}-\varepsilon}}-1\right) 
\geq  \delta n \frac{y^{\frac{1}{2}+\varepsilon}}{\log y}-\frac{y}{\log y}. \]
For $y=n^2$ this contradicts equation (\ref{product-argument}),
proving that for most primes $p < n^2$ the powers $1,2,4,\ldots , 2^n$
 occupy more than
$ y^{\frac{1}{2}-\varepsilon}>  p^{\frac{1}{2}-\varepsilon}$ distinct residue classes
modulo $p$. 
%Hence $\ord_p(2)$ must be at least $p^{1/2-\varepsilon}$.
For more general discussions on this see
Elsholtz \cite{Elsholtz:2002}.

We now come to the additive ingredients, which are more modern, and much 
deeper.

\begin{lemma}[Hart {\cite{Hart:2013}}]{\label{lem:Hart1}}

 Let $R \subseteq \F_p^*$ be a multiplicative subgroup
 such that $|R| \ge p^{\kappa}$, where $\kappa > \frac{11}{23}$.
 Then for all sufficiently large $p$ we have 
 $6R:=R+R+R+R+R+R\supseteq \F_p^{*}$.
 \end{lemma}

An earlier version with exponent $41/83$ appeared in
Schoen and Shkredov {\cite[Theorem 4.1]{SchoenandShkredov:2011}},
and  $55/112$ in Shkredov {\cite[Corollary 32]{Shkredov:2012}}.

See also
\cite{CochraneandPinner:2008, CochraneHartPinnerSpencer:2014, Elsholtz:2008}
for a number of related results on Waring's problem modulo $p$.

\begin{lemma}[Hart \cite{Hart:2013}, Theorem 10]{\label{lem:Hart}}
Let $R$ be a multiplicative subgroup of  $\F_p^*$ with $|R| \ll 
p^{5/9- \varepsilon}$, then $|R+R| \gg |R|^{8/5} (\log |R|)^{-3/10}$. 
\end{lemma}

\begin{lemma}[Schoen and Shkredov 
{\cite[Theorem 2.6, 
$l=3$]{SchoenandShkredov:2011}}]{\label{lemma:Schoenshkredov3}}

Let $R$ be a multiplicative subgroup of $\F_p^*$, with $|R|> p^{2/3}$.
Then $3R:=R+R+R\supseteq \F_p^*$.
\end{lemma}
\begin{lemma}[Schoen and Shkredov 
{\cite[Theorem 2.6, $l=2$]{SchoenandShkredov:2011}}]{
\label{lemma:Schoenshkredov2}}
Let $R$ be a multiplicative subgroup of $\F_p^*$, with $|R|> p^{3/4}$.
Then $2R:=R+R\supseteq \F_p^*$.
\end{lemma}

For the case of a twofold sum there was important earlier work by
Heath-Brown and Konyagin \cite{Heath-BrownandKonyagin:2000}:
Let $R$ be a multiplicative subgroup of $\F_p^*$, with $|R|\gg p^{2/3}$.
Then $|R+R|\gg p$.
Also, several related results of this type of theorem are
due to Bourgain \cite{Bourgain:2005},
Glibichuk \cite{Glibichuk:2008}, \cite{Glibichuk:2011}, 
Cochrane and Pinner \cite{CochraneandPinner:2008}, 
Hart and Iosevich \cite{HartandIosevich:2008}.

\subsection{Proofs of Theorems}
After this preparation
the Theorems are straightforward consequences:\\
The sets $R=\{r^i:1 \leq i \leq \ord_p(r) \}$,
$R=\{r^is^j:1 \leq i \leq \ord_p(r),
1 \leq j \leq \ord_p(s) \}$ etc. are multiplicative subgroups of $\F_p^*$ so
that the corresponding lemmas apply as follows:

Theorem \ref{thm:6fold}
follows from Lemma \ref{lem:erdosmurty}, with $d=1$, $b_1=r$ and
Lemma \ref{lem:Hart1}. 
(Observe that $\frac{11}{23}<\frac{1}{2}$.)
Further observe that it suffices to estimate how often
$\ord_p(r) \leq p^{11/23+ \varepsilon}$ holds. It
 follows from Lemma \ref{lem:Pappalardi}
with $d=1, b_1=r, t=x^{11/23 +\varepsilon/2}$ that
\[ \# \{ p: \ord_p(r) \leq p^{11/23 +\varepsilon/2}\}
\leq  \# \{ p: \ord_p(r) \leq x^{11/23 +\varepsilon/2} \}
\ll_r  \frac{x^{22/23 +\varepsilon}}{\log x}. \]

Theorem \ref{thm:hart-2fold}
follows from Lemma \ref{lem:erdosmurty} with $d=1$, $b_1=r$ and 
Lemma \ref{lem:Hart} with $|R|\geq p^{1/2+ \varepsilon}$, and observing that
$\frac{p^{8\varepsilon/5}}{(\log p)^{3/10}}\gg p^{\varepsilon}$.

\noindent Theorem \ref{thm:3fold}
follows from Lemma \ref{lem:erdosmurty}, with $d=2$, and
Lemma \ref{lemma:Schoenshkredov3}.\\
Theorem \ref{thm:2fold}
follows from Lemma \ref{lem:erdosmurty}, with $d=3$, and
Lemma \ref{lemma:Schoenshkredov2}.\\
Corollary  \ref{thm:2foldconditional}
follows from Lemma \ref{lemma:Schoenshkredov2}.

\subsection{Proof of Theorem {\ref{thm:Skalba2-extension}}}

Every multiple of an integer of the form $M=2^n-1$ has Hamming weight
at least $n$, see Stolarsky {\cite[Theorem 2.1]{Stolarsky:1980}}
or Wagstaff \cite{Wagstaff:2001}.
 Moreover, if 
\begin{equation}{\label{eq:mersennefactors}}
R:=\Omega(2^n-1)< \frac{\log n}{\log k},
\end{equation}
 then some prime factor $q\mid 2^n-1$ also has this property.
 To see this we adapt Ska\l ba's argument.
Let $M=2^n-1=\prod_{i=1}^r q_i^{t_i}$. 
Suppose all $q_i$ have some multiple $w_iq_i$ which
is a sum of at most $k$ powers of two, that is 
\[ w_i q_i = 
\varepsilon_{1,i} 2^{s_1(q_i)}+ \cdots +\varepsilon_{k,i}2^{s_k(q_i)} ,\]
with $\varepsilon_{j,i}\in \{0,1\}$, (not all of them being $0$), then
\[ P:=\prod_{i=1}^r (w_i q_i )^{t_i} = \prod_{i=1}^r 
\left( \varepsilon_{1,i} 2^{s_1(q_i)} + \cdots 
+\varepsilon_{k,i}2^{s_k(q_i)} \right)^{t_i}\]
is on the one hand side a multiple of $M$, 
and thus has Hamming weight at least $n$. 
On the other side, $P$ is a sum 
of at most $k^{\sum_{i=1}^r t_i}=k^{\Omega(2^n-1)}=k^R<n$ powers of 2, 
a contradiction.
Hence there is some prime divisor $q\mid 2^n-1$ for which
\[
\varepsilon_{i} 2^{s_1(q)}+ \cdots +\varepsilon_{k}2^{s_k(q)}\not\equiv 0 \bmod
q,\]
and hence all nonzero multiples of $q$ have Hamming weight at least $k+1$.

\subsection{Proof of theorem \ref{thm:rudnevcounterexample}}
\begin{proof}
In the proof below the $c_i$ are suitably chosen positive constants.

By Heath-Brown's work 
\cite[Corollary 3]{Heath-Brown:1986}
on Artin's primitive root conjecture, 
for any fixed bound $x$ 
and any four distinct positive 
square-free numbers $a_1,a_2,a_3,a_4$ (not being 1) 
(or for any three distinct primes)
one of these, say $a$, is a primitive root (i.e.~$\ord_p(a)=p-1$)
for at least $c_1 \frac{x}{(\log x)^2}$ such primes $p \leq x$.
Hence for suitable $x$ 
there are at least $c_2 \frac{x}{(\log x)^2}$ such primes 
$\frac{x}{2}< p \leq x$. 
Write $p-1=2^s w$, where $w$ is odd, and 
$y=\frac{c_2}{4} \frac{p}{(\log p)^3}$. 
As there are
at most $y \frac{\log x}{\log 2} < \frac{c_2}{2} \frac{x}{(\log x)^2}$ 
such primes with $w\leq y$, 
there are at least $\frac{c_2}{2}\frac{x}{(\log x)^2}$ 
many primes $p$ with $\ord_p(a^{2^s})=w\gg \frac{p}{(\log p)^3}$.
Observe that all elements of the form $a^{2^s\, i}$ have odd order and
that therefore $-1$ is not of this type.
As $a^{2^s\, i}+1\not\equiv 0 \bmod p$ it follows:
the set $A=\{a^{2^s\, i}: 1 \leq i\leq \ord_p(a^{2^s})\}$ is of
 size $|A|\gg \frac{p}{(\log p)^3}$ but
$0 \not\in A+A=A^2+A^2$.

\end{proof}

\section{Open problems}
Finally, we state as open questions:\\
\begin{enumerate}
\item
For primes $p \not\in T_2$: does
 $\ord_p(2) \leq c_{\varepsilon} p^{{1/2}+ \varepsilon}$ hold, for all
$\varepsilon >0$, and some constant $c_{\varepsilon}$? 
%(The table contains a column of the proportion
%$\frac{\ord_p(2)}{\sqrt{p}}$.)
\item
For $h \rightarrow \infty$, what can one say about 
primes $p \not\in T_h$, but not being Mersenne numbers?
\item How can one algorithmically find a ``sparse'' representation
of a multiple of $p$?
Can this be used for other complexity questions?
\end{enumerate}

Finally, we note that there are 231 primes $p \leq 4 \cdot \, 10^6$ 
which are not in
$T_2$, whereas there are $283\, 146$ primes below 4 million.
On the other hand, 
as conjectured by Ska\l ba \cite{Skalba:2004}, proving that there are
 infinitely many such primes 
 remains an open problem.

\section{Acknowledgements}
I would like to thank Rainer Dietmann, Igor Shparlinski and the referees 
for useful comments.

\ \\
\textrm{
%MSC2010:\\ 
\noindent 11A41 Primes\\
11B13   	Additive bases, including sumsets\\
11B50   	Sequences (mod $m$)\\
11N25   	Distribution of integers with specified multiplicative
constraints\\
11P05   	Waring's problem and variants}

\end{document}